\documentclass[11pt,reqno]{amsart}
\usepackage{graphicx}

\usepackage{amscd,amssymb,amsmath,amsthm}
\usepackage{graphicx}
\usepackage{color}
\usepackage{cite}
\newtheorem{thm}{Theorem}
\newtheorem{defn}{Definition}
\newtheorem{lemma}{Lemma}

\numberwithin{equation}{section} 

\def\R{\mathbb{R}}

\begin{document}

\title[On a two dimensional dynamical system]{On a two dimensional dynamical system generated by the floor function}

\author{Usmonov J.B.}

\address{J.\ B.\ Usmonov \\ Institute of mathematics,
81, M.Ulugbek str., 100125, Tashkent,\, \, \, \, \, \, Uzbekistan.}
\email {javohir0107@mail.ru}

\begin{abstract}
In this paper we investigate the two dimensional dynamical system generated by the floor function with a parameter $\lambda\in \R$. 
We describe all limit points of the dynamical system depending on $\lambda$ and on the initial point.
\end{abstract}

\maketitle

\section{Introduction}

Let $X\subset \mathbb R^{2}$ and $A$ be an operator from $X$ to itself. The main problem of the theory of dynamical systems is to study iterations of $A$ : if $A^{n}$ denotes the $n$-fold composition (iteration) of $A$ with itself, then for a given point $x$ one investigates the sequence $x, Ax, A^{2}x, A^{3}x$, and so on. This sequence is called two-dimensional discrete time dynamical system or the forward orbit of $x$, or just the orbit of $x$ for short ([1], [3]).

\begin{defn}
A point $z\in \mathbb R^{2}$ is called a fixed point of $A$ operator if $A(z)=z$. The set of all fixed points is denoted by ${\rm Fix}(A)$.
\end{defn}

\begin{defn}
 The point $z\in \mathbb R^{2}$ is a periodic point of period $n$ if $A^{n}(z)=z$ and $A^{n-1}(z)\neq z$. We denote the set of periodic points of period $n$ by ${\rm Per}_{n}(A)$. The set of all iterates of a periodic point form a periodic orbit.
\end{defn}

For a given operator $A:\mathbb R^{2}\rightarrow\mathbb R^{2}$ the $\omega$-limit set of $z\in\mathbb R^{2}$, denoted by $\omega(z,A)$ or $\omega(z)$, is the set of limit points of the forward orbit $\{A^{n}(z)\}_{n\in\mathbb N}$. Hence, $y\in \omega(z)$ if and only if there is a strictly increasing sequence of natural numbers $\{n_{k}\}_{k\in\mathbb N}$ such that $A^{n_{k}}(z)\rightarrow y$ as $k\rightarrow\infty$.

In this paper we will study the dynamical system generated by operator
$$A: z=(x,y)\in \mathbb R^{2}\rightarrow z'=(x',y')=A(z)\in \mathbb R^{2}$$
defined by
 $$A(z)=\left\{
        \begin{array}{ll}
          x'=\lfloor \lambda y\rfloor & \hbox{} \\
          y'=\lfloor\lambda x\rfloor & \hbox{}
        ,\end{array}
      \right.
  z=(x,y)\in \mathbb R^{2}$$
where $\lambda \in \mathbb R$  is parameter and $\lfloor x\rfloor$ denotes the integer part of $x$.

In our case the dynamical system is
$$z=(x,y), \ \ A(z)=(\lfloor\lambda y\rfloor,\lfloor\lambda x\rfloor), \ \ A^{2}(z)=(\lfloor\lambda\lfloor\lambda x\rfloor\rfloor,\lfloor\lambda\lfloor\lambda y\rfloor\rfloor), \ \ \ldots$$
The main problem is to investigate the following limit
$$\lim_{n\rightarrow\infty}A^{n}(z),$$
for any $z\in \mathbb R^{2}$.
\section{The main results}
\subsection[Fixed points]{Fixed points}

The following lemma gives all fixed points of this operator.

\begin{lemma}
 \emph{For the set of fixed points the followings hold:
 \begin{itemize}
  \item[1)] If $\lambda <0 (\lambda \neq -1)$, then ${\rm Fix}(A)=\{(0, 0)\}$;
  \item[2)] If $\lambda =-1$,  then ${\rm Fix}(A)=\{(m, -m)\, |\, m \in \mathbb Z\} $;
  \item[3)] If $\frac{m-1}{m}<\lambda \leq \frac{m}{m+1}$  for some $m \in \mathbb N$, then 
  $${\rm Fix}(A)=\{(x, \lfloor\lambda x\rfloor)\, |\, x\in \{0,-1,-2,...,-m\}\};$$
  \item[4)] If $\lambda =1$, then ${\rm Fix}(A)=\{(m,m)\,|\,m \in \mathbb Z\}$;
  \item[5)] If $\frac{m+1}{m}\leq \lambda < \frac{m}{m-1}$ for some $m \in \mathbb N$, then 
  $${\rm Fix}(A)=\{(x, \lfloor\lambda x\rfloor)\,|\,x\in \{0,1,2,...,m-1\}\}.$$
\end{itemize}}
\end{lemma}
 {\sf Proof.} For finding fixed points of the operator we need solve $A(z)=z$, i.e., that system $\left\{\begin{array}{ll}
        x=\lfloor\lambda y\rfloor & \\
        y=\lfloor\lambda x\rfloor.&\end{array}\right.$
 By expressing $y$ we get the equation $\lfloor\lambda\lfloor\lambda x\rfloor\rfloor=x$. Roots of $\lfloor\lambda x\rfloor=x$ are also roots of $\lfloor\lambda\lfloor\lambda x\rfloor\rfloor=x$. That's why fixed points given in parts 1-5 of Lemma took by solving $\lfloor\lambda x\rfloor=x$ (see e.g. [2]). But there may be some fixed points of $g(x)=\lfloor\lambda\lfloor\lambda x\rfloor\rfloor$ those are roots of $\lfloor\lambda \lfloor\lambda x\rfloor\rfloor=x$ and are not roots of $\lfloor\lambda x\rfloor=x$. We shall prove that fixed points of operator are just roots of  $\lfloor\lambda x\rfloor=x$.

 In [2] all limit points $\omega(x)$ of the floor function $f(x)=\lfloor\lambda x\rfloor$ found for $\forall \lambda, x\in R$. It was proved that the limit of $\{f^{n}(x)\}$ converges: to fixed points; to $\infty$, or -$\infty$(case $\lambda=-1$ exception). For the case $\lambda=-1$ it was proved that ${\rm Per}_{2}(f)=\mathbb Z$.

  By ${\rm Per}_{2}(f)\subset {\rm Fix}(g)$ we have part 2 of Lemma. Because of $\{g^{n}(x)\}=\{f^{2n}(x)\}$, if $\{f^{n}(x)\}$ is convergent then $\{g^{n}(x)\}$ also is convergent, if $\{f^{n}(x)\}$ converges to $\infty$ then $\{g^{n}(x)\}$ also converges to $\infty$. Thus fixed points of $g(x)$ for $\lambda\neq-1$ consist of roots of $\lfloor\lambda x\rfloor=x$ only. Now we prove part 2, $x=\lfloor-y\rfloor=\lfloor-\lfloor-x\rfloor\rfloor=x$.

\subsection{The limit points}

Now we shall describe the set $\omega(z)$ for each given $z\in\mathbb R^{2}$.
\subsubsection{\bf The case $\lambda\leq 0$.}
\begin{thm}
If $\lambda <0$, then the dynamical system generated by operator $A$ has the following properties:
\begin{itemize}
  \item[(1)] If $-1<\lambda <0$, then 
  $$\lim_{n\rightarrow \infty}A^{n}(z)=(0,0)$$ 
  for all $z \in\mathbb R^{2}$.
  \item[(2)] If $\lambda =-1$, then each pairs of integer numbers has period two and
$$\omega(z)=\left\{\begin{array}{ll}
\{z,A(z)\}, \ \ \ \ \ \ \ \mbox{if} \ \ z\in \mathbb Z\times\mathbb Z; \\
\{A(z),A^{2}(z)\}, \ \ \mbox{if} \ \ z\in\mathbb R^{2}\setminus\mathbb  Z\times\mathbb Z.
\end{array}\right.$$
 \item[(3)] If $\lambda <-1$,  then $A(z)=0$ for $z\in U_{\frac{1}{|\lambda|}}^{-}(0)=\{(x,y) \, | \, \frac{1}{\lambda}<x\leq0, \frac{1}{\lambda}<y\leq0\}$ and
$$\omega(z)=\left\{\begin{array}{ll}
\{(0,0)\}, \ \ \ \ \ \ \ \ \ \ \ \ \ \ \ if \ \ z\in U^{-}_\frac{1}{|\lambda|}(0); \\[3mm]
\{(\infty,\infty),(-\infty,-\infty)\}, \ \  if \ \ z\in\{\mathbb R_{++}^{2}\cup\mathbb R_{--}^{2}\}\setminus U^{-}_\frac{1}{|\lambda|}(0);\\[3mm]
\{(\infty,-\infty),(-\infty,\infty)\}, \ \ if \ \ z\in\{\mathbb R_{+-}^{2}\cup\mathbb R_{-+}^{2}\}\setminus U_{\frac{1}{|\lambda|}}^{-}(0).
\end{array}\right.$$
\end{itemize}
\end{thm}
 where \\
$\mathbb R_{++}^{2}=\{(x,y)\, |\,x,y\in\mathbb R,x>0,y>0\}$,\ \ $\mathbb R_{--}^{2}=\{(x,y)\,\ |\,\ x,y\in\mathbb R,x<0,y<0\}$,\\[2mm]
$\mathbb R_{+-}^{2}=\{(x,y)\,\ |\,\ x,y\in\mathbb R,x>0,y<0\}$,\ \ $\mathbb R_{-+}^{2}=\{(x,y)\,\ |\,\ x,y\in\mathbb R,x<0,y>0\}$.
\\[2mm]
{\sf Proof.} (1) Let $-1<\lambda <0$. For all $z=(x_{0},y_{0})\in\mathbb R^{2}$ we have
\begin{equation}\label{1}
z=(x_{0},y_{0}),\ \ A(z)=(\lfloor\lambda y_{0}\rfloor,\lfloor\lambda x_{0}\rfloor),\ \ A^{2}(z)=(\lfloor\lambda\lfloor\lambda x_{0}\rfloor\rfloor,\lfloor\lambda\lfloor\lambda y_{0}\rfloor\rfloor),\ \ ... .\end{equation}
We can separate sequence (\ref{1}) to two subsequences $\{u_{n}\}$ and $\{v_{n}\}$:
\begin{equation}\label{2}
u_{n}=\left\{\begin{array}{ll}
x_{n}, \ \  \hbox{if} \ \ n-\hbox{even} \\
y_{n}, \ \  \hbox{if} \ \ n-\hbox{odd}
\end{array}\right.
\end{equation}
\begin{equation}\label{3}
v_{n}=\left\{\begin{array}{ll}
x_{n}, \ \ \hbox{if} \ \ n-\hbox{odd} \\
y_{n}, \ \ \hbox{if} \ \ n-\hbox{even}.
\end{array}\right.
\end{equation}

Let's write several terms of those sequences,
$$\{u_{n}\}: u_{0}=x_{0}, \ \ u_{1}=y_{1}=\lfloor\lambda x_{0}\rfloor, \ \ u_{2}=x_{2}=\lfloor\lambda\lfloor\lambda x_{0}\rfloor\rfloor, \ \ u_{3}=y_{3}=\lfloor\lambda\lfloor\lambda\lfloor\lambda x_{0}\rfloor\rfloor\rfloor,\ \ ...$$
$$\{v_{n}\}: v_{0}=y_{0}, \ \ v_{1}=x_{1}=\lfloor\lambda y_{0}\rfloor, \ \ v_{2}=y_{2}=\lfloor\lambda\lfloor\lambda y_{0}\rfloor\rfloor, \ \ v_{3}=x_{3}=\lfloor\lambda\lfloor\lambda\lfloor\lambda y_{0}\rfloor\rfloor\rfloor,\ \ ...$$

 In Theorem 2 of [2] was proved $\lim_{n\rightarrow\infty}u_{n}=0$ for any $\lambda\in(-1,0)$ and for all $x_0 \in\mathbb R$.
 That's why we have
$$\lim_{n\rightarrow \infty}A^{n}(z)=(0,0).$$

(2) If $\lambda=-1$, then $A^{2}(z)=z$ and $A(z)\neq z$ for all $z\in\mathbb Z\times\mathbb Z$. Thus each pairs of integer numbers has period two.
If $z \in \mathbb R^{2}\setminus Z\times\mathbb Z$, then $A(z)\in \mathbb Z\times\mathbb Z$. So we have
$$\omega(z)=\left\{\begin{array}{ll}
\{z,A(z)\}, \ \ \ \ \ \ \ \mbox{if} \ \ z\in \mathbb Z\times\mathbb Z; \\[3mm]
\{A(z),A^{2}(z)\}, \ \ \mbox{if} \ \ z\in\mathbb R^{2}\setminus\mathbb  Z\times\mathbb Z.
\end{array}\right.$$

(3) Let $\lambda<-1$ and $z=(x_0,y_0) \in U^{-}_{\frac{1}{|\lambda|}}(0)$. Then we have $\lfloor\lambda x_{0}\rfloor=0$ and $\lfloor\lambda y_{0}\rfloor=0$, i.e., $A(z)=0$. If $x_{0}\leq \frac{1}{\lambda}$ then for sequence (\ref{2}) we have
$$u_{1}<|u_{2}|\leq u_{3}<|u_{4}|\leq u_{5}<\ldots$$
for all $z \in\mathbb R^{2}\setminus U^{-}_{\frac{1}{|\lambda|}}(0) $ and, if $x_{0}>0$ then
$$|u_{1}|\leq u_{2}<|u_{3}|\leq u_{4}<|u_{5}|\leq \ldots .$$
Since $\{|u_{n}|\}_{n\geq 1}\subset\mathbb N$  and $\lim_{n\rightarrow\infty}|u_{n}|=\infty$ the following hold
$$
\lim_{n\rightarrow\infty}A^{2n}(z)=\left\{\begin{array}{ll}
(\infty,\infty), \ \ \ \ \ \ \hbox{if} \ \ z \in\mathbb R^{2}_{+,+}\setminus U^{-}_{\frac{1}{|\lambda|}}(0); \\[3mm]                                     (-\infty,-\infty), \ \ \hbox{if} \ \ z \in\mathbb R^{2}_{-,-}\setminus U^{-}_{\frac{1}{|\lambda|}}(0); \\[3mm]
(\infty,-\infty),\ \ \ \ \hbox{if} \ \ z \in\mathbb R^{2}_{+,-}\setminus U^{-}_{\frac{1}{|\lambda|}}(0); \\[3mm]
(-\infty,\infty),\ \ \ \ \hbox{if} \ \ z \in\mathbb R^{2}_{-,+}\setminus U^{-}_{\frac{1}{|\lambda|}}(0).
\end{array}\right.$$

$$\lim_{n\rightarrow\infty}A^{2n+1}(z)=\left\{\begin{array}{ll}
(-\infty,-\infty),\ \ \  \hbox{if} \ \ z \in\mathbb R^{2}_{+,+}\setminus U^{-}_{\frac{1}{|\lambda|}}(0); \\[3mm]
(\infty,\infty),\ \ \ \ \ \ \  \hbox{if} \ \ z \in\mathbb R^{2}_{-,-}\setminus U^{-}_{\frac{1}{|\lambda|}}(0); \\[3mm]
(-\infty,\infty), \ \ \ \ \hbox{if} \ \ z \in\mathbb R^{2}_{+,-}\setminus U^{-}_{\frac{1}{|\lambda|}}(0); \\[3mm]
(\infty,-\infty),\ \ \ \ \hbox{if} \ \ z \in\mathbb R^{2}_{-,+}\setminus U^{-}_{\frac{1}{|\lambda|}}(0).
\end{array}\right.$$

\subsubsection{\bf The case $0<\lambda<1$.}

Note that for each $\lambda\in (0,1)$ there exists $m\in\mathbb N$ such that ${m-1\over m}<\lambda\leq{m\over m+1}$.

\begin{thm}
    Let $\frac{m-1}{m}<\lambda\leq\frac{m}{m+1}$ for some $m\in\mathbb N$. Then the following hold:
\begin{itemize}
  \item[(1)] If $z \in \{(x_{0},\,y_{0})\,  |\,  x_{0}\geq 0,\, y_{0}\geq 0\}$, then
$$\lim_{n\rightarrow\infty}A^{n}(z)=(0,0).$$
  \item[(2)] If $z \in \{(x_{0},y_{0})\, |\, \frac{k}{\lambda}\leq x_{0}<\frac{k+1}{\lambda},y_{0}\geq0\}\bigcup\{(x_{0},y_{0})\, |\, \frac{k}{\lambda}\leq y_{0}<\frac{k+1}{\lambda},x_{0}\geq0\}$, then
$$\omega(z)=\{(k,0),(0,k)\};$$
where $k\in\{-1,-2,...,-m\}.$
  \item[(3)] If $z \in \{(x_{0},y_{0})\, |\, x_{0}<\frac{-m}{\lambda},y_{0}\geq0\}\bigcup \{(x_{0},y_{0})\, |\, y_{0}<\frac{-m}{\lambda},x_{0}\geq0\}$, then
$$\omega(z)=\{(-m,0),(0,-m)\}.$$
    \item[(4)] If $z \in \{(x_{0},y_{0})\, |\, x_{0}<\frac{-m+1}{\lambda},y_{0}<\frac{-m+1}{\lambda}\}$, then
$$\lim_{n\rightarrow\infty}A^{n}(z)=(-m,-m).$$
  \item[(5)] If $z \in \{(x_{0},y_{0})\, |\, \frac{k}{\lambda}\leq x_{0}<\frac{k+1}{\lambda},\frac{p}{\lambda}\leq y_{0}<\frac{p+1}{\lambda}\}$, then
$$\omega(z)=\{(k,p),(p,k)\};$$
where $k\in\{-1,-2,...,-m\}, p\in\{-1,-2,...,-m\}$.
  \item[(6)] If $z \in \{(x_{0},y_{0})\, |\, x_{0}<\frac{-m+1}{\lambda},\frac{k}{\lambda}\leq y_{0}<\frac{k+1}{\lambda}\}\bigcup \{(x_{0},y_{0})\, |\, \frac{k}{\lambda}\leq x_{0}<\frac{k+1}{\lambda},y_{0}<\frac{-m+1}{\lambda}\}$, then
$$\omega(z)=\{(k,-m),(-m,k)\};$$
where $k\in\{-1,-2,...,-m\}$.
\end{itemize}
\end{thm}
 {\sf Proof.} (1) For subsequences (\ref{2}) and (\ref{3}) of (\ref{1}) we have
$$0\leq \lim_{n\rightarrow \infty}u_{n}\leq\lim_{n\rightarrow\infty}\lambda^{n}u_{0}=0,$$
$$0\leq \lim_{n\rightarrow \infty}v_{n}\leq\lim_{n\rightarrow\infty}\lambda^{n}v_{0}=0,$$
for all $z=(x_{0},y_{0})\in\mathbb R^{2}_{++}$, i.e. $\lim_{n\rightarrow\infty}u_{n}=\lim_{n\rightarrow\infty}v_{n}=0$. Then $\lim_{n\rightarrow\infty}A^{n}(z)=(0,0)$.

(2) If $z\in \{(x_{0},y_{0})\ \ |\ \ \frac{k}{\lambda}\leq y_{0}<\frac{k+1}{\lambda},y_{0}\geq0\}$ then $\lfloor\lambda x_{0}\rfloor=k$, where $k\in\{-1,-2,...,-m\}$. Since $k$ is a fixed point of $f(x)=\lfloor\lambda x\rfloor$, then  $\lim_{n\rightarrow\infty}u_{n}=k$ and  by proof of 1st part we have $\lim_{n\rightarrow\infty}v_{n}=0$ for $y_{0}\geq0$.

In case $\forall z\in \{(x_{0},y_{0})\, |\, \frac{k}{\lambda}\leq y_{0}<\frac{k+1}{\lambda},\, x_{0}\geq0\}$ we can write  $\lim_{n\rightarrow\infty}u_{n}=0$ and $\lim_{n\rightarrow\infty}v_{n}=k$ as above. That means, $\omega(z)=\{(k,0),(0,k)\}$.

(3) $u_{1}<-m$ and $u_{1}>u_{0}$ are true for all $z \in \{(x_{0},y_{0})\, |\, x_{0}<\frac{-m}{\lambda},\, y_{0}\geq0\}$ . For terms of (2) we see that $u_{n+1}>u_{n}$, i.e. $\{u_{n}\}$ is an increasing sequence, which is bounded from above by $-m$. Since $-m$ is the unique fixed point of $f(x)=\lfloor\lambda x\rfloor$ in $(-\infty,-m]$, we have $\lim_{n\rightarrow\infty}u_{n}=-m$. For $\{v_{n}\}$ we have $\lim_{n\rightarrow\infty}v_{n}=0$.

In case $z \in \{(x_{0},y_{0})\, |\, y_{0}<\frac{-m}{\lambda},x_{0}\geq0\}$ we can write $\lim_{n\rightarrow\infty}u_{n}=0$ and $\lim_{n\rightarrow\infty}v_{n}=-m$, then $\omega(z)=\{(-m,0),(0,-m)\}$.

(4) Like the proof of 3rd part  we may write $$\lim_{n\rightarrow\infty}u_{n}=\lim_{n\rightarrow\infty}v_{n}=-m\Rightarrow\lim_{n\rightarrow\infty}A^{n}(z)=(-m,-m),$$
for all $z \in \{(x_{0},y_{0})\, |\, x_{0}<\frac{-m+1}{\lambda},y_{0}<\frac{-m+1}{\lambda}\}$.

(5) If $z\in\{(x_{0},y_{0})\, |\, \frac{k}{\lambda}\leq x_{0}<\frac{k+1}{\lambda}, \frac{p}{\lambda}\leq x_{0}<\frac{p+1}{\lambda}\}$ then $\lfloor\lambda x_{0}\rfloor=k, \lfloor\lambda y_{0}\rfloor=p$. Since $k$ and $p$ are fixed points of $f(x)=\lfloor\lambda x\rfloor$ we have
$$\lim_{n\rightarrow\infty}u_{n}=k, \lim_{n\rightarrow\infty}v_{n}=p\Rightarrow \omega(z)=\{(k,p),(p,k)\},$$
where $k\in\{-1,-2,...,-m\}, p\in\{-1,-2,...,-m\}$.

(6) The proof is based on parts 3-4.

\subsubsection{\bf The case $\lambda\geq1$.}

In case $\lambda=1$, the form of operator $A$ is $A(z)=\left\{
                                                   \begin{array}{ll}
                                                     x'=\lfloor y\rfloor & \hbox{} \\
                                                     y'=\lfloor x\rfloor & \hbox{}
                                                   \end{array}
                                                 \right.
$ and
$$\omega(z)=\{(\lfloor x\rfloor,\lfloor y\rfloor),(\lfloor y\rfloor,\lfloor x\rfloor)\}$$
for all $z=(x,y)\in\mathbb R^{2}$.

 \begin{thm}
 Let $\frac{m+1}{m}\leq\lambda<\frac{m}{m-1}$ for some $m\in\mathbb N$. Then the following hold:
\begin{itemize}
  \item[(1)] If $z\in\{(x_{0},y_{0})\, |\, \frac{k}{\lambda}\leq x_{0}<\frac{k+1}{\lambda},y_{0}<0\}\bigcup \{(x_{0},y_{0})\, |\, \frac{k}{\lambda}\leq y_{0}<\frac{k+1}{\lambda},x_{0}<0\}$, then
$$\omega(z)={(k,-\infty),(-\infty,k)};$$
where $k\in\{0,1,2,...,m-1\}$.
  \item[(2)] If $z\in\{(x_{0},y_{0})\, |\, x_{0}<0,y_{0}\geq\frac{m}{\lambda}\}\bigcup \{(x_{0},y_{0})\, |\, y_{0}<0,x_{0}\geq\frac{m}{\lambda}\}$, then
$$\omega(z)=\{(\infty,-\infty),(-\infty,\infty)\}.$$
  \item[(3)] If $z\in \{(x_{0},y_{0})\, |\, \frac{k}{\lambda}\leq x_{0}<\frac{k+1}{\lambda},\frac{p}{\lambda}\leq x_{0}<\frac{p+1}{\lambda}\}$, then
$$\omega(z)=\{(k,p),(p,k)\};$$
where $k\in\{0,1,2,...,m-1\}, p\in\{0,1,2,...,m-1\}$.
  \item[(4)] If $z\in\{(x_{0},y_{0})\, |\, \frac{k}{\lambda}\leq x_{0}<\frac{k+1}{\lambda},y_{0}\geq\frac{m}{\lambda}\}\bigcup\{(x_{0},y_{0})\, |\, \frac{k}{\lambda}\leq y_{0}<\frac{k+1}{\lambda},x_{0}\geq\frac{m}{\lambda}\}$, then
$$\omega(z)=\{(k,\infty),(\infty,k)\};$$
where $k\in\{0,1,2,...,m-1\}$.
  \item[(5)] If $z\in\{(x_{0},y_{0})\, |\, x_{0}<0,y_{0}<0\}$, then
$$\lim_{n\rightarrow\infty}A^{n}(z)=(-\infty,-\infty).$$
  \item[(6)] If $z\in\{(x_{0},y_{0})\, |\, x_{0}\geq\frac{m}{\lambda},y_{0}\geq\frac{m}{\lambda}\}$, then
$$\lim_{n\rightarrow\infty}A^{n}(z)=(\infty,\infty).$$
\end{itemize}
\end{thm}
 {\sf Proof.} (1) Since $z\in\{(x_{0},y_{0})\, |\, \frac{k}{\lambda}\leq x_{0}<\frac{k+1}{\lambda},y_{0}<0\}$ then $\lim_{n\rightarrow\infty}u_{n}=k$.
$y_{0}>\lambda y_{0}\geq\lfloor\lambda y_{0}\rfloor=f(y_{0})$ because $y_{0}<0$ and $\lambda y_{0}<0$  for $\lambda>1$. Using this inequality we get $f^{n}(y_{0})>f^{n+1}(y_{0}) (y_{n}>y_{n+1})$. Due to lack of fixed points of $f(x)=\lfloor\lambda x\rfloor$ in $(-\infty,0)$, we have $\lim_{n\rightarrow\infty}f^{n}(y_{0})=\lim_{n\rightarrow\infty}v_{n}=-\infty$.

In this case also $\forall z\in\{(x_{0},y_{0})\, |\, \frac{k}{\lambda}\leq y_{0}<\frac{k+1}{\lambda},x_{0}<0\}$, we may write $\lim_{n\rightarrow\infty}v_{n}=k$ and $\lim_{n\rightarrow\infty}u_{n}=-\infty$ as above. Thus $\omega(z)=\{(k,-\infty),(-\infty,k)\}$.

(2) We showed that $\lim_{n\rightarrow\infty}u_{n}=-\infty$ is true for all $z\in\{(x_{0},y_{0})\, |\, x_{0}<0,y_{0}\geq\frac{m}{\lambda}\}$ in part 2. If $y_{0}\geq\frac{m}{\lambda}$, we have $v_{n}<v_{n+1}$ and $\{v_{n}\}$ is an increasing sequence, that's why
$\lim_{n\rightarrow\infty}v_{n}=\infty$. In case when $z\in\{(x_{0},y_{0})\, |\, y_{0}<0,x_{0}\geq\frac{m}{\lambda}\}$ we get $\lim_{n\rightarrow\infty}u_{n}=\infty$ and $\lim_{n\rightarrow\infty}v_{n}=-\infty$, i.e. $\omega(z)=\{(-\infty,\infty),(\infty,-\infty)\}$.

(3) If $z\in \{(x_{0},y_{0})\, |\, \frac{k}{\lambda}\leq x_{0}<\frac{k+1}{\lambda},\frac{p}{\lambda}\leq x_{0}<\frac{p+1}{\lambda}\}$ then $\lfloor\lambda x_{0}\rfloor=k$ and $\lfloor\lambda y_{0}\rfloor=p$. Since  $k$ and $p$ are  fixed points of $f(x)=\lfloor\lambda x\rfloor$ , we get $\lim_{n\rightarrow\infty}u_{n}=k, \lim_{n\rightarrow\infty}v_{n}=p\Rightarrow \omega(z)=\{(k,p),(p,k)\}$.

Proofs of parts 4-6 are directly come from above results.

\end{document}